\documentclass[12,reqno]{amsart}
\usepackage{amsmath,amsthm,amsfonts,amssymb,times}
\usepackage[mathscr]{eucal}
\usepackage[breaklinks]{hyperref}
\usepackage{verbatim}
\usepackage{url}
\usepackage{xcolor}
\setbox0=\hbox{$+$}
\newdimen\plusheight
\plusheight=\ht0
\def\+{\;\lower\plusheight\hbox{$+$}\;}

\setbox0=\hbox{$-$}
\newdimen\minusheight
\minusheight=\ht0
\def\-{\;\lower\minusheight\hbox{$-$}\;}

\setbox0=\hbox{$\cdots$}
\newdimen\cdotsheight
\cdotsheight=\plusheight
\def\cds{\lower\cdotsheight\hbox{$\cdots$}}

\makeatletter
\def\leqalignno#1{\displ@y \tabskip\z@ plus\@ne fil
	\halign to\displaywidth{\hfil$\@lign\displaystyle{##}$\tabskip\z@skip
		&$\@lign\displaystyle{{}##}$\hfil\tabskip\z@ plus\@ne fil
		&\kern-\displaywidth\rlap{$\@lign\hbox{\rm##}$}\tabskip\displaywidth\crcr
		#1\crcr}}
\makeatother

\newcommand{\eb}{\begin{equation}}
	\newcommand{\ee}{\end{equation}}

\newcommand{\df}{\dfrac}

\renewcommand{\Re}{\operatorname{Re}}

\newcommand{\s}{{\sigma}}

\renewcommand{\a}{\alpha}
\renewcommand{\b}{\beta}

\newcommand{\g}{\gamma}

\renewcommand{\l}{\lambda}

\renewcommand{\t}{\tau}

\renewcommand{\Re}{\textup{Re}}

\renewcommand{\(}{\left\(}
\renewcommand{\)}{\right\)}
\renewcommand{\[}{\left\[}
\renewcommand{\]}{\right\]}
\renewcommand{\i}{\infty}

\numberwithin{equation}{section}
\theoremstyle{plain}
\newtheorem{theorem}{Theorem}[section]
\newtheorem{lemma}[theorem]{Lemma}
\newtheorem{corollary}[theorem]{Corollary}

\newtheorem{remark}[theorem]{Remark}

\newtheorem*{theorem*}{Theorem}

\allowdisplaybreaks

\begin{document}
	\title[A segment of Euler product]{A segment of Euler product associated to a certain Dirichlet series} 
	\author{Rajat Gupta}
	\address{Institute of Mathematics, Academia Sinica, 6F, Astronomy-Mathematics Building, No. 1, Sec. 4, Roosevelt Road, Taipei 106319, TAIWAN\\ \newline
	\textbf{Current address:} The university of Texas at Tyler, Texas 75701, USA.}
	\email{rajatgpt972@gmail.com, rgupta@uttyler.edu}
	\author{Aditi Savalia}
	
	\thanks{2010 \textit{Mathematics Subject Classification.} Primary 11M41; Secondary, 11M06.\\
		\textit{Keywords and phrases.} Euler product,  Dirichlet series, approximate functional equation, and divisor function}
	\address{Department of Mathematics, Indian Institute of Technology Bombay, Mumbai 400076, Maharashtra, India} 
	\email[Corresponding author]{aditiben.s@iitgn.ac.in, aditisav1995@gmail.com }
	
	\maketitle
	
	\begin{abstract}
		In the spirit of the work of Hardy-Littlewood and Lavrik, we study the Dirichlet series associated to the generalized divisor function $\sigma_{\a}(n):=\sum_{d|n}d^{\alpha}$. We obtain an exact identity relating the Dirichlet series $\zeta(s)\zeta(s-\alpha)$ and a segment of the Euler product attached to it. Specifically, our main theorems are valid in the critical strip.
	\end{abstract}
	
	\section{Introduction}
	
	
	In the early $20$th century, Hardy and Littlewood \cite{gl3} studied the approximation of the Riemann zeta function $\zeta(s):=\sum_{n=1}^{\infty}1/{n^s}$ in the region $0\leq$ $\textup{Re}(s)\leq 1$. They expressed it as a sum of two finite series with two error terms, called an approximate functional equation,  namely for $0 \leq \sigma \leq 1$, 
	\begin{align}\label{app1}
		\zeta(s)=\sum_{n\leq x}\frac{1}{n^s}+\chi(s)\sum_{n\leq y}\frac{1}{n^{1-s}}+\mathcal{O}\left(\frac{1}{x^\sigma} \right)+\mathcal{O}\left(\frac{|t|^{\frac{1}{2} -\sigma}}{y^{1-\sigma}} \right). 
	\end{align}
	Here, $\chi(s):=2(2\pi)^{s-1}\sin\left(\pi s/2\right)\Gamma\left(1-s\right)$, $s=\sigma+it$ with $|t|=2\pi x y,$ and $x,y> A$ for any $A>0$. They called this identity a ``compromise" between the series expression of $\zeta(s)$ and its functional equation $\zeta(s)=\chi(s)\zeta(1-s)$. In the same paper they also derived the approximation for $\zeta^2(s)$, 
	namely that 
	\begin{align}\label{app2}
		\zeta^2(s)=\sum_{n\leq x}\frac{d(n)}{n^s}+\chi^2(s)\sum_{n< y}\frac{d(n)}{n^{1-s}}+\mathcal{O}\left(x^{1/2-\sigma}\log t \right),
	\end{align}
for $0 \leq \sigma \leq 1$,	where $x,y>A$ and $xy=\left(t/2\pi\right)^2$. Note that error terms in both \eqref{app1} and \eqref{app2} are best possible (see \cite[Chapter 4]{iv}) in the sense that for some values of $x,y,$ and $t$, the corresponding error terms are of the exact order of magnitude. Later, Heath-Brown \cite[Lemma 1]{hb1}, \cite[Lemma 1]{hb2}, and Lavrik \cite{lavrik-1} studied the weighted approximate functional equations. For example, they derived them for any positive integer power of $\zeta(s)$ similar to what we have in \eqref{app1} and \eqref{app2}. In \cite{Simonic}, Simoni\v{c} calculated implied constants in the error terms in \eqref{app1} and \eqref{app2}.

	Lavrik \cite{lavrik5} made further investigations to approximate $\zeta(s)$. Interestingly, rather than truncating the series representation of $\zeta(s)$, he considered the truncation of its Euler product. Indeed, he  studied a segment of the Euler product which, in turn, extends the Euler product formula 
	\begin{align}
		\label{Ep}
		\zeta(s)=\prod_{p}\left(1-\frac{1}{p^{s}} \right)^{-1}
	\end{align}
	into the critical strip. 
	We record his result \cite[Equation (1)]{lavrik5} here. Throughout the paper we are fixing the notation that $\mathbb{P}$ and $\bar{\mathbb{P}}$ denote a finite set of primes and their product, respectively.
	
	\begin{theorem}[{Lavrik, \cite[Equation (1)]{lavrik5}}] Let $\mathbb{P}\neq\emptyset$, $\Gamma(s)$ be the gamma function and denote $(n,\bar{\mathbb{P}})$ to be the greatest common divisor of $n$ and $\bar{\mathbb{P}}$. Then for complex numbers $s$ and $\tau$, where \textup{Re}$(s)>0$ \footnote{Here we have corrected the condition from any complex $s$ to \text{Re}$(s)>0$ which Lavrik missed while evaluating integration of $n=1$ term coming from $(6)$ in \cite{lavrik5}.} with $s \neq 1$ and \textup{Re}$(\tau)>0$, we have
	\begin{align} \label{eq: Maintheorem1eqn2}
		\zeta(s)\prod_{p \in \mathbb{P}}\left(1-\frac{1}{p^s} \right)&=1+\frac{(\pi\tau)^{s/2}}{\Gamma\left( s/2\right)}\Bigg(\frac{\t^{-1/2}}{s-1}\prod_{p \in \mathbb{P}}\left(1-\frac{1}{p} \right)  \nonumber\\
		&+\int_{1}^{\i}t^{s/2-1}\Bigg(\sum\limits_{\mathop {n=2 {}}\limits_{(n,\bar{\mathbb{P}})=1} }^{\i}e^{-\pi\tau n^2 t} +\frac{1}{t^s}\Bigg(\sum\limits_{\mathop {n=2 {}}\limits_{(n,\bar{\mathbb{P}})=1} }^{\i}e^{-\pi\tau n^2/t}-\frac{1}{2}\sqrt{\frac{t}{\t}}\prod_{p \in \mathbb{P}}\left(1-\frac{1}{p} \right)\Bigg)\Bigg)dt\Bigg).
	\end{align} 
\end{theorem}
	As an immediate application of \eqref{eq: Maintheorem1eqn2}, in a follow-up paper
	, he gave an approximation \cite[Theorem 1]{lavrik6} of $\zeta(s)$ in terms of a finite segment of the Euler product and then used it to give a new method of deriving a zero free region for $\zeta(s)$. We refer readers to \cite[Section 1.2]{lavrik6} for further discussion on the connection between the approximate functional equation and the non-trivial zeros of the Riemann zeta function. Approximate functional equations are one of the very effective tools in analytic number theory through which we can express any $L$-function, in a region where its Dirichlet series does not converge, by a summation of  finite Dirichlet polynomials up to some error term.

	Our main goal in this paper is to express a Dirichlet series associated to a generalized divisor function $\sigma_{\a}(n):=\sum_{d|n}d^{\alpha}$ in terms of a segment of an Euler product in the critical strip where the Dirichlet series and its Euler product are not convergent, see our main result given in Theorem \ref{Maintheorem1}. Indeed, we derive as a special case of Theorem \ref{Maintheorem1} at $\alpha \to 0$, the following identity for $\zeta^{2}(s)$ in terms of a segment of an Euler product and a series involving the modified Bessel function of the second kind of order zero. We record the identity here. 
	\begin{theorem}\label{Theoalpha0eqnintro}
		Let $\mathbb{P} \neq \emptyset$, $\textup{Re}(\t)>0$, and \textup{Re}$(s)>0$ with $s \neq 1$. Then
		\begin{align}\label{coroalpha0eqnintro}
			&\zeta^2(s)\prod_{p \in \mathbb{P}}\left(1-\frac{1}{p^s} \right)^2 \nonumber \\
			&=1+\frac{\pi^{s}\tau^{s/2}}{\Gamma^2\left( \frac{s}{2}\right)}\Bigg[\frac{1}{(s-1)\sqrt{\t}}\prod_{p \in \mathbb{P}}\left(1-\frac{1}{p} \right)\left(-2\sum_{n|\bar{\mathbb{P}}}\frac{\mu(n)}{n}\log n +\left(\frac{1}{s-1}+\gamma -\log(4\pi \sqrt{\t})\right)\prod_{p \in\mathbb{P}}\left(1-\frac{1}{p} \right)\right)\nonumber\\
			&+2\int_{1}^{\i}t^{s/2-1}\Bigg(\sum\limits_{\mathop {n=2 {}}\limits_{(n,\bar{\mathbb{P}})=1} }^{\i}d(n)K_{0}\left( 2n\pi \sqrt{\t t}\right) +t^{-s}\Bigg(\sum\limits_{\mathop {n=2 {}}\limits_{(n,\bar{\mathbb{P}})=1} }^{\i}d(n)K_{0}\left( 2n\pi\sqrt{\frac{\t}{t}}\right)\nonumber\\
			&-\frac{1}{4}\sqrt{\frac{t}{\t}}\prod_{p \in \mathbb{P}}\left(1-\frac{1}{p} \right)\Bigg(-2\sum_{n|\bar{\mathbb{P}}}\frac{\mu(n)}{n}\log n +\left(\gamma-\log\left(4\pi \sqrt{\frac{\t}{t}}\right)\right)\prod_{p \in\mathbb{P}}\left(1-\frac{1}{p} \right)\Bigg)\Bigg)\Bigg)dt\Bigg].
		\end{align}
	\end{theorem}
	
	We adapt Lavrik's approach to achieve our goal.
	In particular, we derive a formula of type \eqref{eq: Maintheorem1eqn2} for $\zeta(s)\zeta(s-\a)$ and its special cases, one of which is given in Theorem \ref{Theoalpha0eqnintro}. To structure the proof of our main theorem in this paper, we will be employing Ramanujan's  renowned version of the Guinand formula  \cite[p.~253]{berndtleesohn}, which he had discovered several years before Guinand \cite{apg2} himself.
	The formula reads as follows:
	
	\begin{theorem}[{B.C.~Berndt, Y.~Lee, and J.~Sohn \cite[p.~253]{berndtleesohn}}]
	Let $ K_{\nu}(z) $ denote the modified Bessel function of the second kind of order $ \nu $, and $ x $, $ y $ are positive numbers such that $ xy=\pi^{2} $.  For $ z $ being any complex number, we have
	\begin{align}\label{apgeqn}
		&\sqrt{x}\sum_{n=1}^{\i}\s_{-z}(n)n^{z/2}K_{z/2}(2nx)-\sqrt{y}\sum_{n=1}^{\i}\s_{-z}(n)n^{z/2}K_{z/2}(2ny) \nonumber\\
		&\qquad \qquad \qquad=\frac{1}{4}\Gamma\left( \frac{z}{2}\right)\zeta(z)\left(y^{(1-z)/2}-x^{(1-z)/2}\right)+\frac{1}{4}\Gamma\left( -\frac{z}{2}\right)\zeta(-z)\left(y^{(1+z)/2}-x^{(1+z)/2}\right).
	\end{align}
	\end{theorem}

	\section{Preliminaries}
	
	We refer readers to G.~N.~Watson's classical treatise for the definition and properties of the modified Bessel function of the second kind $K_{\nu}(z)$, see \cite[pp.~15, 78]{watson}.  The following lemmas will be used in the sequel.

	\begin{lemma}\label{asymptotic}\cite[pp.~202]{watson} Let $K_{\nu}(x)$ denote the modified Bessel function of order $\nu$. As $z\to\infty$,
		\begin{align}
			K_{\nu}(z)=&\left(\df{2}{\pi z}\right)^{1/2}e^{-z}\left(1+O\left(\df{1}{z}\right)\right).\label{asymptotic2}
		\end{align}
	\end{lemma}
	
	
	Next, the Mellin transform \cite{kp} of a function is a multiplicative analog of the Fourier transform. The Mellin transform of a locally integrable function $f (x)$ on $(0,\infty)$ is given by 
	\begin{align}\label{mellindef}
		F(s)=\int_{0}^\infty f(x)x^{s-1}~dx,
	\end{align} when the integral converges. Let the region of convergence of the above integral \eqref{mellindef} be $a<\Re (s)=c <b$. Then the inversion formula for Mellin transformation gives 
	\begin{align*}
		f(x)=\frac{1}{2\pi i}\int_{c-i\infty }^{c+i\infty}F(s)x^{-s}~ds, \qquad (a<c<b).
	\end{align*} 
	We will be using a well-known result called Parseval's formula \cite[Equation (3.1.13)]{kp} noted below.
	\begin{lemma}\label{Parsevalle's idty} Let $F(s)$ and $G(s)$ be the Mellin transform of $f$ and $g$ respectively. Then for $c$ lying in the common strip of analyticity of $F$ and $G$,
		\begin{align*}\int_{0}^\infty f(u/x) g(x)~\frac{dx}{x} =\frac{1}{2\pi i}\int_{c-i\infty}^{c+i\infty}F(s)G(s)u^{-s}~ds.\end{align*}
	\end{lemma}
	
	We conclude this section by stating three lemmas which are quite useful to prove our main theorems.  
	\begin{lemma}\label{Mellin 1.0}
		For $ c>\max\{0,-\textup{Re}(s),-\textup{Re}(s)+\textup{Re}(\a)\} $ and $\textup{Re}(u^2)>0$, we have
		\begin{align*}
			\frac{1}{2\pi i } \int_{c -i \i}^{c+i \i}\Gamma\left( \frac{z+s}{2}\right)\Gamma\left(\frac{z+s-\a}{2}\right)u^{-z}\dfrac{dz}{z}=
			2u^{s-\a/2}\int_{1}^{\i}t^{-\a/4+s/2-1}K_{\a/2}\left( 2u\sqrt{t}\right)dt.
		\end{align*}
	\end{lemma}
	We prove the above Lemma \ref{Mellin 1.0} in  Section \ref{4}.
	
	 { Let the Dirichlet convolution of two arithmetic functions $f,g$ be defined 
		as
		\begin{align*}
			f\ast g(n)=\sum_{d|n}f(d)g(n/d).
		\end{align*}
		Also, we consider $\mu$ to be the well-known M\"obius function defined as $\mu(1)=1$ and
		\begin{align*}
			\mu(n)= \left\{\begin{array}{cc}
				(-1)^k & \text{if $n$ is a product of $k$ distinct primes, }\\
				0 & \text{otherwise.}
			\end{array}\right.
		\end{align*} 
	Moreover, we let $I$ be the multiplicative identity function, namely $I(1) = 1$ and $I(n) = 0$ for $n > 1$.
	We also remind the reader that $(\mu N^\alpha)(n)= \mu(n)n^\alpha$ and that $\sigma_{\a}^{-1}(n)= (\mu \ast \mu N^\alpha)(n)$ in the sense of the Dirichlet inverse (refer \cite[Theorem $2.20$]{apostol}).
	We state the following lemma without its proof as it trivially follows using elementary methods of number theory.}

	
	\begin{lemma}\label{identity for mu}
		For $z, \a \in \mathbb{C},$ we have
		\begin{align*}
			\sum_{k|\bar{\mathbb{P}}^{2}}\frac{(\mu \ast \mu N^\a)~(k)}{k^{z}}=\prod_{p \in \mathbb{P}}\left(1-\frac{1}{p^{z}} \right) \left(1-\frac{1}{p^{z-\a}} \right).
		\end{align*}
	\end{lemma}

	\begin{lemma}\label{Galpha} For $\mathbb{P}\neq \emptyset$ and $(n,\bar{\mathbb{P}}^2)\neq 1$,
		\begin{align}\label{eq: g_alpha}
			\mathcal{G}_\a(n):=\sum_{k|(n,\bar{\mathbb{P}}^2)}\s_\a \left( \frac{n}{k}\right)\left( \mu \ast \mu N^\a\right)(k)=0.
		\end{align} 
	\end{lemma}

\begin{proof} 
	This can be seen as follows. Let $\sigma_{\a}^{-1}$ denote the Dirichlet inverse of $\sigma_\a$, i.e., $\sigma_{\a}\ast\sigma_{\a}^{-1}=I$. Let $n$ be such that $p^3|n$ for some prime $p$. Then
	\begin{align*}
		\s_{\a}^{-1}(n)=(\mu \ast \mu N^\a)~(n)=\sum_{d|n}\mu\left(\frac{n}{d}\right)\mu\left(d\right)d^\a=0,
	\end{align*}
	as at least one of the $d$ or $n/d$ is not square free.
	Thus $\s_\a^{-1}(n)=0$ for such $n$. In the other case, i.e., when there is no $p$ with $p^3|n$, we can write $n=\bar{n}\cdot t$, such that $p|\bar{n} \Rightarrow p\in\mathbb{P}$, and $(t,\bar{\mathbb{P}})=1$, indeed, $	\bar{n}=(n,\bar{\mathbb{P}}^2)$. Moreover we are in the case where $(n,\bar{\mathbb{P}}^2)>1$, thus we have 
	\begin{align*}
		\sum_{k|(n,\bar{\mathbb{P}}^2)}\s_\a \left( \frac{n}{k}\right)\s_{\a}^{-1}(k)=\s_\a(t)\sum_{k|\bar{n}}\s_\a \left( \frac{\bar{n}}{k}\right)\s_{\a}^{-1}(k)=\s_\a(t)\big(\s_\a * \s_{\a}^{-1}\big)(\bar{n})=0.
	\end{align*}
\end{proof}



	\section{Main Theorems and Intermediate results}\label{3}
	
	We begin this section by recording an interesting modular type relation. Recall the notation that $\mathbb{P}$ and $\bar{\mathbb{P}}$ denote a finite collection of primes and their product, respectively.
	\begin{theorem}\label{doublesumBessel}
		For $\mathbb{P} \neq \emptyset$, $\textup{Re}(t/\t)>0$, and $\a \in \mathbb{C}$, we have
		\begin{align}\label{123}
			&\sum_{k|\bar{\mathbb{P}}^2}\frac{(\mu \ast \mu N^\a)~(k)}{k^{1+\a/2}}\sum_{n=1}^{\i}\frac{\s_\a(n)}{n^{\a/2}}K_{\a/2}\left( \frac{2n\pi}{k}\sqrt{\frac{t}{\tau}}\right)\nonumber\\
			&=\sqrt{\frac{\t}{t}}\sum\limits_{\mathop {n=1 {}}\limits_{(n,\bar{\mathbb{P}})=1} }^{\i}\frac{\s_\a(n)}{n^{\a/2}}K_{\a/2}\left( 2n\pi\sqrt{\frac{\t}{t}}\right) -\left(\mathcal{F}_{\a}(t)+\mathcal{F}_{-\a}(t)\right),
		\end{align}
		where,
		\begin{align}\label{FG}
			\mathcal{F}_{\a}(t)= \frac{1}{4}\Gamma\left( \frac{\a}{2}\right)\zeta(\a)\left(\frac{1}{\pi} \sqrt{\frac{\t}{t}}\right)^{\a/2}\prod_{p \in \mathbb{P}}\left(1-\frac{1}{p} \right)\left(1-\frac{1}{p^{1-\a}} \right).
		\end{align}
	\end{theorem}
	Using the above modular relation we will derive our main theorem stated below.
	\begin{theorem}\label{Maintheorem1}
		For $\mathbb{P} \neq \emptyset$, $\textup{Re}(\t)>0$, $\a \in \mathbb{C}\backslash \{0\}$, and \textup{Re}$(s)>$ $\max\{0,\textup{Re}(\alpha)\}$ with $s\neq 1, 1+\a$, we have
		\begin{align}\label{Maintheorem1eqn}
			&\zeta(s)\zeta(s-\a)\prod_{p \in \mathbb{P}}\left(1-\frac{1}{p^s} \right)\left(1-\frac{1}{p^{s-\a}} \right)=1+\frac{\pi^{s-\a/2}\tau^{s/2}}{\Gamma\left( \frac{s}{2}\right)\Gamma\left( \frac{s-\a}{2}\right)}\Bigg(\frac{\t^{-1/2}}{s-1}\mathcal{R}_{1-\a} +\frac{\t^{-(1+\a)/2}}{s-1-\a}\mathcal{R}_{1+\a} \nonumber \\
			&\qquad \qquad +2\tau^{-\a/4}\int_{1}^{\i}t^{-\a/4+s/2-1}\Bigg(\sum\limits_{\mathop {n=2 {}}\limits_{(n,\bar{\mathbb{P}})=1} }^{\i}\frac{\s_\a(n)}{n^{\a/2}}K_{\a/2}\left( 2n\pi \sqrt{\t t}\right) \nonumber
			\\
			& \qquad \qquad \qquad\qquad+\frac{t^{\a/2}}{t^s}\Bigg(\sum\limits_{\mathop {n=2 {}}\limits_{(n,\bar{\mathbb{P}})=1} }^{\i}\frac{\s_\a(n)}{n^{\a/2}}K_{\a/2}\left( 2n\pi\sqrt{\frac{\t}{t}}\right)
			-\sqrt{\frac{t}{\t}}\left(\mathcal{F}_\a(t)+\mathcal{F}_{-\a}(t) \right)\Bigg)\Bigg)dt\Bigg),
		\end{align} 
		where
		\begin{align*}
			\mathcal{R}_\a =\pi^{-\a/2}\Gamma\left(\frac{\a}{2}\right)\zeta(\a)\prod_{p \in \mathbb{P}}\left(1-\frac{1}{p} \right)\left(1-\frac{1}{p^{\a}} \right),
		\end{align*}
		and $\mathcal{F}_{\a}(t)$ is defined in \eqref{FG}.
	\end{theorem}
	It can be easily seen from Lemma \ref{asymptotic} that both the series in the integrand on the right-hand side are absolutely convergent. Also the absolute convergence of the integral on the right-hand side above can be justified using the modular relation given in  Theorem \ref{doublesumBessel}.  Indeed, with the help of \eqref{123}, the integrand can be written as multiplication of the factors $t^{-\a/4+s/2-1}$ and
	\begin{align*}
		\sum\limits_{\mathop {n=2 {}}\limits_{(n,\bar{\mathbb{P}})=1} }^{\i}\frac{\s_\a(n)}{n^{\a/2}}K_{\a/2}\left( 2n\pi \sqrt{\t t}\right)+\frac{t^{\alpha/2-s+1/2}}{\sqrt{\tau}}\sum_{k|\bar{\mathbb{P}}^2}\frac{(\mu \ast \mu N^\a)~(k)}{k^{1+\a/2}}\sum_{n=1}^{\i}\frac{\s_\a(n)}{n^{\a/2}}K_{\a/2}\left( \frac{2n\pi}{k}\sqrt{\frac{t}{\tau}}\right).
	\end{align*}
	Therefore, again using Lemma \ref{asymptotic}, we can show that the integrand has an exponential decay. 
	
	We will show that Theorem \ref{Theoalpha0eqnintro} is a limiting case $\a \to 0$ of Theorem \ref{Maintheorem1}. Further restricting ourselves to the region $\textup{Re}(s)>1$, Theorem \ref{Theoalpha0eqnintro} reduces to:
	\begin{corollary}\label{coroalpha0>1}
		For $\mathbb{P} \neq \emptyset$, $\textup{Re}(s)>1$, and $\textup{Re}(\t) >0,$ we obtain
		\begin{align}\label{eq}
			\zeta^2(s)\prod_{p \in \mathbb{P}}\left(1-\frac{1}{p^s} \right)^2=1&+\frac{2\pi^{s}\tau^{s/2}}{\Gamma^2\left( \frac{s}{2}\right)}\int_{1}^{\i}t^{s/2-1}\Bigg(\sum\limits_{\mathop {n=2 {}}\limits_{(n,\bar{\mathbb{P}})=1} }^{\i}d(n)K_{0}\left( 2n\pi \sqrt{\t t}\right)\nonumber\\
			&\qquad+t^{-s}\sum\limits_{\mathop {n=2 {}}\limits_{(n,\bar{\mathbb{P}})=1} }^{\i}d(n)K_{0}\left( 2n\pi\sqrt{\frac{\t}{t}}\right)\Bigg)dt. 
		\end{align}

	\end{corollary}
	
	\begin{remark}
		Using \cite[Formula~6.561.16, p.~676]{grn}, we can further simplify the integration in \eqref{eq} and derive the following elementary identity for $\textup{Re}(s)>1$,
			\begin{align}
			\sum\limits_{\mathop {n=1 {}}\limits_{(n,\bar{\mathbb{P}})=1} }^{\i}\frac{d(n)}{n^s}=\zeta^2(s)\prod_{p \in \mathbb{P}}\left(1-\frac{1}{p^s} \right)^2.
		\end{align}
		If we let $\mathbb{P}$ be the collection of all primes in the above corollary, then we recover the well-known Euler identity, namely, we get for $\textup{Re}(s) > 1$,
		$$\zeta^2(s)\prod_{p}\left(1-\frac{1}{p^s}\right)^2=1.$$
	\end{remark}
	

	\section{Proofs of Lemmas, main theorem and their corollaries}\label{4}
	\subsection{Proof of Lemma \textup{\ref{Mellin 1.0}} }
	\begin{proof}[Proof of Lemma \textup{\ref{Mellin 1.0}}]
		For $\textup{Re}(z)>0$, we have \cite[p.~347, Equation (1)]{ober}
		\begin{align}\label{GammaMellin}
			\int_{0}^{\i}x^{z-1}e^{-x}dx= \Gamma(z).
		\end{align}
		Replacing $z$ with $\frac{z+s}{2}$, and then applying change of variable $x \to x^2$ in \eqref{GammaMellin}, we get
		\begin{align}
			2\int_{0}^{\i}x^{s+z-1}e^{-x^2}dx=\Gamma\left(\frac{z+s}{2} \right),~\textup{Re}(z)>-\textup{Re}(s). \label{Mellin 1}
		\end{align}
		Also, for $\textup{Re}(\l+\b)>0$, $\textup{Re}(\mu)>0$, and $\textup{Re}(\mu+\nu)>0$, we have \cite[p.~657,~Formula 6.455(1)]{grn} 
		\begin{align}\label{Mellin 2}
			\int_{0}^{\i}x^{\mu-1}e^{-\b x}\Gamma(\nu, \l x)dx=\frac{\l^{\nu}\Gamma(\mu+\nu)}{\mu(\l+\b)^{\mu+\nu}}~{}_{2}F_{1}\left(1,\mu+\nu;\mu+1;\frac{\b}{\l+\b} \right).
		\end{align}
		Letting $\b \to 0$, $x \to x^2$, $\mu =z/2,~\l =1$, and $\nu =\frac{s-\a}{2}$ in \eqref{Mellin 2}, we obtain for $\textup{Re}(z)>\max\{0, -\textup{Re}(s)+\textup{Re}(\a)\}$
		\begin{align}
			\int_{0}^{\i}x^{z-1}\Gamma\left(\frac{s-\a}{2}, x^2\right)dx=\frac{1}{z}\Gamma\left(\frac{z+s-\a}{2} \right),\label{Mellin 1.1}
		\end{align}
		where $ \Gamma\left( a,x\right)  $ is the incomplete gamma function defined as 
		\begin{align*}
			\Gamma\left( a,x\right) := \int_{x}^{\infty} t^{a-1}e^{-t}dt.
		\end{align*}
		Employing Parseval's formula stated in Lemma \ref{Parsevalle's idty},  
		for $c$ in the intersection of $\textup{Re}(z)>0$, $\textup{Re}(z)>-\textup{Re}(s)$ and $\textup{Re}(z)>-\textup{Re}(s)+\textup{Re}(\a)$, the identities \eqref{Mellin 1} and \eqref{Mellin 1.1} gives
		\begin{align*}
			\frac{1}{2\pi i}\int_{c-i\infty}^{c+i\infty}\Gamma\left(\frac{z+s}{2} \right)\Gamma\left(\frac{z+s-\a}{2} \right)\frac{u^{-z}}{z}dz&= 2\int_{0}^{\i}\Gamma\left(\frac{s-\a}{2}, \frac
			{u^2}{x^2}\right)e^{-x^2}x^{s}\frac{dx}{x} \nonumber \\
			&=2\int_{0}^{\i}\left(\int_{u^2/x^2}^{\i}t^{\frac{s-\a}{2}-1}e^{-t}dt\right)e^{-x^2}x^{s-1}~dx.
		\end{align*}
		We make the change of variable $t \to \frac{u^2}{x^2}t$, to obtain
		\begin{align}
			\frac{1}{2\pi i}\int_{c-i\infty}^{c+i\infty}\Gamma\left(\frac{z+s}{2} \right)\Gamma\left(\frac{z+s-\a}{2} \right)\frac{u^{-z}}{z}dz&=2\int_{0}^{\i}\left(\int_{1}^{\i}\left(\frac{u^2}{x^2} t\right)^{\frac{s-\a}{2}-1}e^{-u^2 t/x^2}dt\right)\frac{u^2}{x^2}e^{-x^2}x^{s-1}~dx\nonumber\\
			&=2u^{s-\a}\int_{0}^{\i}\left(\int_{1}^{\i}t^{\frac{s-\a}{2}-1}e^{-u^2 t/x^2}dt\right)e^{-x^2}x^{\a-1}~dx \nonumber\\
			&=2u^{s-\a}\int_{1}^{\i}\left(\int_{0}^{\i}e^{-u^2 t/x^2-x^2}x^{\a-1}~dx\right)t^{\frac{s-\a}{2}-1}dt,\label{changeofintehration}
		\end{align}
		where in the last step, change of order of integration can be justified by  invoking Watson's lemma \cite[Theorem 2.4]{temme}. 
		We have from \cite[p.~368,~3.471,~Formula (9)]{grn},
		\begin{align*}
			\int_{0}^{\i}x^{\nu -1}e^{-\gamma x -\b/x}dx=2\left( \frac{\b}{\g}\right)^{\nu/2}K_{\nu}\left(2\g\sqrt{\b} \right),  
		\end{align*}
		for $ \textup{Re}(\b)>0 $ and $ \textup{Re}(\g)>0 $.	
		We now use the above integral with change of variable $x \to x^2$ followed by taking $\b=u^2t$, $\gamma =1$, and $\nu= \frac{\a}{2}$ in \eqref{changeofintehration} to derive our lemma, namely, for $ \textup{Re}(z)>\max\{0,-\textup{Re}(s),-\textup{Re}(s)+\textup{Re}(\a)\} $ and $\textup{Re}(u^2)>0$,
		\begin{align*}\
			\frac{1}{2\pi i}\int_{c-i\infty}^{c+i\infty}\Gamma\left(\frac{z+s}{2} \right)\Gamma\left(\frac{z+s-\a}{2} \right)\frac{u^{-z}}{z}dz=2u^{s-\a/2}\int_{1}^{\i}t^{-\a/4+s/2-1}K_{\a/2}\left(2u\sqrt{t} \right)dt.
		\end{align*}
	\end{proof}

	\subsection{Proofs of Theorems \textup{\ref{doublesumBessel}} and \textup{\ref{Maintheorem1}}}
	\begin{proof}[Proof of Theorem \textup{\ref{doublesumBessel}}]
		
		In \textit{Ramanujan-Guinand's} formula take $z =\a$, $\displaystyle{x= \frac{\pi}{k\sqrt{(\tau/t)}}}$, and $y =\displaystyle{\frac{\pi^2}{x} =\pi k\sqrt{(\tau/t)}}$, and use the relation $\sigma_{-\a}(n)=n^{-\a}\sigma_{\a}(n)$, to get 
		\begin{align*}
			&\left(\frac{\pi}{k}\sqrt{\frac{t}{\t}} \right)^{1/2}\sum_{n=1}^{\i}\frac{\s_{\a}(n)}{n^{\a/2}}K_{\a/2}\left(\frac{2n\pi}{k}\sqrt{\frac{t}{\t}}\right)-\left(\pi k\sqrt{\frac{\t}{t}} \right)^{1/2}\sum_{n=1}^{\i}\frac{\s_{\a}(n)}{n^{\a/2}}K_{\a/2}\left( 2\pi nk\sqrt{\frac{\t}{t}}\right)\\
			&\qquad \qquad \qquad=\frac{1}{4}\Gamma\left( \frac{\a}{2}\right)\zeta(\a)\left(\left( \pi k\sqrt{\frac{\t}{t}}\right)^{(1-\a)/2}-\left(\frac{\pi}{k}\sqrt{\frac{t}{\t}}\right)^{(1-\a)/2}\right)\\
			& \qquad \qquad \qquad \qquad \qquad+\frac{1}{4}\Gamma\left( -\frac{\a}{2}\right)\zeta(-\a)\left(\left( \pi k\sqrt{\frac{\t}{t}}\right)^{(1+\a)/2}-\left(\frac{\pi}{k}\sqrt{\frac{t}{\t}}\right)^{(1+\a)/2}\right).
		\end{align*}
		After considerable simplification we obtain
		\begin{align}\label{Rg}
			&\sum_{n=1}^{\i}\frac{\s_{\a}(n)}{n^{\a/2}}K_{\a/2}\left(\frac{2n\pi}{k}\sqrt{\frac{t}{\t}}\right)=k\sqrt{\frac{\t}{t}}\sum_{n=1}^{\i}\frac{\s_{\a}(n)}{n^{\a/2}}K_{\a/2}\left( 2\pi nk\sqrt{\frac{\t}{t}}\right)-A_\a -A_{-\a},
		\end{align}
		where
		\begin{align*}
			A_\a=A_\a(t):=-\frac{1}{4}\Gamma\left( \frac{\a}{2}\right)\zeta(\a)\left( k\sqrt{\frac{\t}{t}}\left( \pi k\sqrt{\frac{\t}{t}}\right)^{-\a/2}-\left(\frac{k}{\pi}\sqrt{\frac{\t}{t}}\right)^{\a/2}\right).
		\end{align*}
		We multiply both sides by $\displaystyle{\frac{(\mu \ast \mu N^\a)~(k)}{k^{1+\a/2}}}$ and then take the sum over divisors of $\bar{\mathbb{P}}^2$. Thus, \eqref{Rg} changes to 
		\begin{align}\label{applyingsumonk}
			\sum_{k|\bar{\mathbb{P}}^2}\frac{(\mu \ast \mu N^\a)~(k)}{k^{1+\a/2}}\sum_{n=1}^{\i}\frac{\s_{\a}(n)}{n^{\a/2}}K_{\a/2}\left(\frac{2n\pi}{k}\sqrt{\frac{t}{\t}}\right)=\sqrt{\frac{\t}{t}}&\sum_{n=1}^{\i}\frac{K_{\a/2}\left( 2n\pi\sqrt{\frac{\t}{t}}\right)}{n^{\a/2}}\sum_{k|(\bar{\mathbb{P}}^2,n)}(\mu\ast \mu N^\a)~(k)\s_\a\left(\frac{n}{k}\right)\nonumber\\
			&-\sum_{k|\bar{\mathbb{P}}^2}\frac{(\mu \ast \mu N^\a)~(k)}{k^{1+\a/2}}\left(A_\a+A_{-\a} \right).
		\end{align}
		The first sum on the right-hand side of the above equation can be split into two sums, namely, 
		\begin{align}\label{split}
			\sum\limits_{\mathop {n=1 {}}\limits_{(n,\bar{\mathbb{P}}^2)=1} }^{\i}\frac{\s_\a(n)}{n^{\a/2}}K_{\a/2}\left( 2n\pi\sqrt{\frac{\t}{t}}\right)
			+\sum\limits_{\mathop {n=1 {}}\limits_{(n,\bar{\mathbb{P}}^2)\neq 1} }^{\i}\frac{\mathcal{G}_\a(n)}{n^{\a/2}}K_{\a/2}\left( 2n\pi\sqrt{\frac{\t}{t}}\right),
		\end{align}
		where $ \mathcal{G}_\a $ is defined in \eqref{eq: g_alpha}. 
		Surprisingly, by employing Lemma \ref{Galpha}, we see that the second sum in \eqref{split} amounts to zero. 
		
		Also note that 
		\begin{align}\label{Aalpha1}
			\sum_{k|\bar{\mathbb{P}}^2}\frac{(\mu \ast \mu N^\a)~(k)}{k^{1+\a/2}}A_\a=&-\frac{1}{4}\Gamma\left( \frac{\a}{2}\right)\zeta(\a)\Bigg(\sqrt{\frac{\t}{\t}}\sum_{k|\bar{\mathbb{P}}^2}\frac{(\mu \ast \mu N^\a)~(k)}{k^{\a}}\left( \pi \sqrt{\frac{\t}{t}}\right)^{-\a/2}\nonumber\\
			&\hspace{4cm}-\sum_{k|\bar{\mathbb{P}^2}}\frac{(\mu \ast \mu N^\a)~(k)}{k}\left(\frac{1}{\pi}\sqrt{\frac{\t}{t}}\right)^{\a/2}\Bigg) \nonumber\\
			=&-\frac{1}{4}\Gamma\left( \frac{\a}{2}\right)\zeta(\a)\Bigg(\pi^{-\a/2}\left(\sqrt{\frac{\t}{t}}\right)^{1-\a/2}\sum_{n|\bar{\mathbb{P}}}\frac{\mu(n)}{n^\a}\sum_{m|\bar{\mathbb{P}}}\mu(m)\nonumber\\
			&\hspace{3.5cm}-\left(\frac{1}{\pi}\sqrt{\frac{\t}{t}}\right)^{\a/2}\sum_{n|\bar{\mathbb{P}}}\frac{\mu(n)}{n}\sum_{m|\bar{\mathbb{P}}}\frac{\mu(m)m^\a}{m}\Bigg) \nonumber\\
			=&\frac{1}{4}\Gamma\left( \frac{\a}{2}\right)\zeta(\a)\left(\frac{1}{\pi}\sqrt{\frac{\t}{t}}\right)^{\a/2}\prod_{p \in \mathbb{P}}\left(1-\frac{1}{p} \right)\left(1-\frac{1}{p^{1-\a}} \right)=\mathcal{F}(\a),
		\end{align}
		where in the last step we use that for $\mathbb{P} \neq \emptyset$, 
		\begin{align*}
			\sum_{m|\bar{\mathbb{P}}}\mu(m) =0,~\text{ and }~ \sum_{n|\bar{\mathbb{P}}}\frac{\mu(n)}{n}\sum_{m|\bar{\mathbb{P}}}\frac{\mu(m)m^\a}{m}=\prod_{p\in\mathbb{P}}\left(1-\frac{1}{p}\right)\left(1-\frac{1}{p^{1-\a}}\right).
		\end{align*}
		Similarly we obtain,
		\begin{align}\label{Aalpha2}
			\sum_{k|\bar{\mathbb{P}}^2}\frac{(\mu \ast \mu N^\a)~(k)}{k^{1+\a/2}}A_{-\a}=\frac{1}{4}\Gamma\left(- \frac{\a}{2}\right)\zeta(-\a)\left(\frac{1}{\pi}\sqrt{\frac{\t}{t}}\right)^{-\a/2}\prod_{p \in \mathbb{P}}\left(1-\frac{1}{p} \right)\left(1-\frac{1}{p^{1+\a}} \right)=\mathcal{F}(-\a)
		\end{align} 
		using the facts that for $\mathbb{P} \neq \emptyset$, 
		\begin{align*}
			\sum_{m|\bar{\mathbb{P}}}\mu(m) =0,~\text{ and }~ \sum_{n|\bar{\mathbb{P}}}\frac{\mu(n)}{n^{1+\a}}\sum_{m|\bar{\mathbb{P}}}\frac{\mu(m)}{m}=\prod_{p\in\mathbb{P}}\left(1-\frac{1}{p}\right)\left(1-\frac{1}{p^{1+
					\a}}\right).
		\end{align*}
		Hence, by plugging \eqref{split}, \eqref{Aalpha1}, and \eqref{Aalpha2} into \eqref{applyingsumonk}, we derive our result.
	\end{proof}
	\begin{proof}[Proof of Theorem \textup{\ref{Maintheorem1}}]
		For $ s=\sigma+it $, and $ \alpha= \sigma_{0}+it_{0} $, let us define $F(z)$ as
		\begin{align*}
			F(z):=\pi^{-z+\alpha/2}\Gamma\left( \frac{z}{2}\right) \Gamma\left( \frac{z-\alpha}{2}\right)\zeta(z)\zeta(z-\alpha)\frac{\tau^{(s-z)/2}}{z-s}\prod_{p \in \mathbb{P}}\left(1-\frac{1}{p^{z}} \right)\left(1-\frac{1}{p^{z-\a}} \right).
		\end{align*}
		
		Let $ \mathcal{C} $ be the contour determined by the line segments $[\sigma_{1}-iT,\sigma_{1}+iT]$ , $[ \sigma_{1}+iT,\sigma_{2}+iT]$, $[\sigma_{2}+iT,\sigma_{2}-iT ]$, $[\sigma_{2}-iT, \sigma_{1}-iT] $, where $ \sigma_{1} < \min\{0,\sigma, \sigma_{0}\} $, $ \sigma_{2}>\max\{1,\sigma,1+\sigma_{0}\} $ and $ T> |t|+|t_0| $.
		
		Note that poles of $ \Gamma\left( \frac{z}{2}\right)$ and $ \Gamma\left( \frac{z-\alpha}{2}\right) $ at $ 0$ and  $\alpha $ are canceled by the zeros of the product over finite collection of primes, and the poles  at $ -2m $ and $ \alpha-2m $ for $ m\in \mathbb{N} $ are canceled by the trivial zeros of $ \zeta(z) $ and $\zeta(z-\alpha) $. For $\textup{Re}(\t) >0$, 
		the factor $\displaystyle\frac{\t^{(s-1)/2}}{z-s}$ also contributes a simple pole at $z= s$. 
		Thus we have poles of $ F(z) $ only at $ z= 1$, $1+\a  $ and $ s $. The choice of $\s_1$, $\s_2$, and $T$ is such that all poles lie inside the contour. Let $R_a$ denote residue of $F(z)$ at the pole $z=a$. Then we have, 
		\begin{enumerate}
			\item at $ z=1 $, $ F(z) $ has a simple pole with residue
			\begin{align*}
				R_{1}&= \lim_{z\rightarrow1}(z-1)F(z)=\pi^{(\alpha-1)/2} \Gamma\left( \frac{1-\alpha}{2}\right)\zeta(1-\alpha)\frac{\tau^{(s-1)/2}}{1-s}\prod_{p \in \mathbb{P}}\left(1-\frac{1}{p} \right)\left(1-\frac{1}{p^{1-\a}} \right),
			\end{align*}
			\item at $ z=1+\a $, $ F(z) $ has a simple pole with residue
			\begin{align*}
				R_{1+a}&= \lim_{z\rightarrow1+\a}(z-1-\a)F(z)=\pi^{-( 1+\a)/ 2}\Gamma\left( \frac{1+\a}{2}\right)\zeta(1+\a)\frac{\tau^{( s-1-\a) /2}}{1+\a-s}\prod_{p \in \mathbb{P}}\left(1-\frac{1}{p} \right)\left(1-\frac{1}{p^{1+\a}} \right),
			\end{align*} 
			\item at $ z=s $, $ F(z) $ has a simple pole with residue
			\begin{align*}
				R_{s}&= \lim_{z\rightarrow s}(z-s)F(z)=\pi^{-s+\a/ 2}\Gamma\left( \frac{s}{2}\right)\Gamma\left( \frac{s-\a}{2}\right)\zeta(s)\zeta(s-\a)\prod_{p \in \mathbb{P}}\left(1-\frac{1}{p^{s}} \right)\left(1-\frac{1}{p^{s-\a}} \right).
			\end{align*} 
		\end{enumerate}
		
		Now applying Cauchy's residue theorem for $F(z)$ in the contour $\mathcal{C}$, we have
		\begin{align}
			\frac{1}{2\pi i }\left(\int_{\s_2 -i T}^{ \s_2 +i T}+\int_{\s_2 +i T}^{ \s_1 +i T}+\int_{\s_1 +i T}^{ \s_1 -i T}+\int_{\s_1 -i T}^{ \s_2 -i T}\right)F(z)dz =R_1 +R_{1+\a}+R_s. \label{sum of residues}
		\end{align}
		
		Using Stirling bound for $\Gamma\left(\frac{z}{2} \right)$, and $\Gamma\left(\frac{z-\a}{2} \right)$ with $|\textup{arg} (\t)|<\pi/2$, which states that for $z'=\sigma'+it'$ with $a\leq \sigma' \leq b$ and $|t'|>1$,
		\begin{align}
			|\Gamma(z')|\ll |t'|^{\sigma'-\frac{1}{2}}e^{-\frac{\pi}{2}|t'|},
		\end{align} 
		and property of the Riemann Zeta function that $\zeta(z')\ll |t'|^{\frac{3}{2}+\delta}$ for $\sigma'>-\delta$, the integrals along the horizontal lines  $[ \sigma_{1}+iT,\sigma_{2}+iT]$ and  $[ \sigma_{1}-iT,\sigma_{2}-iT]$ tend to zero as $T \to \i$. 
		Hence, letting $T \to \i$ in \eqref{sum of residues}, we get 
		\begin{align}\label{cauchy}
			\frac{1}{2\pi i }\left(\int_{\s_2 -i \i}^{ \s_2 +i \i}+\int_{\s_1 +i \i}^{ \s_1 -i \i}\right)F(z)dz =R_1 +R_{1+\a}+R_s.
		\end{align}
		We first evaluate 
		\begin{align}\label{f_1}
			f_{1}:=\frac{1}{2\pi i}\int_{\s_2-i\i}^{\s_2 +i\i}\pi^{-z+\alpha/2}\Gamma\left( \frac{z}{2}\right) \Gamma\left( \frac{z-\alpha}{2}\right)\zeta(z)\zeta(z-\alpha)\frac{\tau^{(s-z)/2}}{z-s}\prod_{p \in \mathbb{P}}\left(1-\frac{1}{p^{z}} \right)\left(1-\frac{1}{p^{z-\a}} \right)dz.
		\end{align}
		By our choice of $ \sigma_{2} $, $ \zeta(z) $ and $ \zeta(z-\a) $ have the Euler product representations on the line segment $[\sigma_{2}+iT,\sigma_{2}-iT ]$. Thus the factor 
		\begin{align}\label{Z(z,t)def}
			Z(z,\a):=\zeta(z)\zeta(z-\a)\prod_{p \in \mathbb{P}}\left(1-\frac{1}{p^{z}} \right)\left(1-\frac{1}{p^{z-\a}} \right), 
		\end{align}
		for $\textup{Re}(z)=\s_2$, reduces to 
		\begin{align}\label{Z(z,t)}
			Z(z,\a)=\sum\limits_{\mathop {n=1 {}}\limits_{(n,\bar{\mathbb{P}})=1} }^{\i}\frac{\s_\a(n)}{n^{z}}.
		\end{align}
		Upon using $\eqref{Z(z,t)}$ in $\eqref{f_1}$, and then changing the order of summation and integration we obtain,
		\begin{align*}
			f_{1}=\frac{1}{2\pi i}\sum\limits_{\mathop {n=1 {}}\limits_{(n,\bar{\mathbb{P}})=1} }^{\i}\frac{\sigma_{\a}(n)}{n^{z}}\int_{\s_2-i\i}^{\s_2 +i\i}\pi^{-z+\alpha/2}\Gamma\left( \frac{z}{2}\right) \Gamma\left( \frac{z-\alpha}{2}\right)\frac{\tau^{(s-z)/2}}{z-s}dz.
		\end{align*}
		%
		By change of variable $ z $ by $ z+s $, we derive, 
		\begin{align}
			f_1 =\pi^{-s+\a/2}\sum_{(n,\bar{\mathbb{P}})=1}\frac{\sigma_{\a}(n)}{n^{s}}\left( \frac{1}{2\pi i } \int_{a -i \i}^{ a+i \i}\Gamma\left( \frac{z+s}{2}\right)\Gamma\left(\frac{z+s-\a}{2}\right)(n\pi\sqrt{\tau})^{-z}\dfrac{dz}{z}   \right), \label{f_1midway}
		\end{align}
		for $a =\s_2 -\s >0$. It is easy to verify that for $\textup{Re}(z)=a$, $\s_2>\max\{1,\s,1+\s_0 \}$ and  $\textup{Re}(\t)>0$, all the conditions are satisfied to employ Lemma \ref{Mellin 1.0} in \eqref{f_1midway} with $ u=n\pi\sqrt{\tau}$. This gives 
	%
	\begin{align*}
		f_1= 2\tau^{s/2-\a/4}\sum\limits_{\mathop {n=1 {}}\limits_{(n,\bar{\mathbb{P}})=1} }^{\i}\frac{\s_\a(n)}{n^{\a/2}}\int_{1}^{\i}t^{-\a/4+s/2-1}K_{\a/2}\left( 2n\pi \sqrt{\t t}\right)dt.
	\end{align*}
	The asymptotic bound on $K_\nu(z)$ from Lemma \ref{asymptotic} allows us to interchange the order of summation and integration, and derive
	\begin{align}\label{f_1main}
		f_1= 2\tau^{s/2-\a/4}\int_{1}^{\i}t^{-\a/4+s/2-1}\sum\limits_{\mathop {n=1 {}}\limits_{(n,\bar{\mathbb{P}})=1} }^{\i}\frac{\s_\a(n)}{n^{\a/2}}K_{\a/2}\left( 2n\pi \sqrt{\t t}\right)dt.
	\end{align}
	Now we turn our attention towards, 
	\begin{align}\label{f_2}
		f_{2}:=\frac{1}{2\pi i}\int_{\s_1+i\i}^{\s_1 -i\i}\pi^{-z+\alpha/2}\Gamma\left( \frac{z}{2}\right) \Gamma\left( \frac{z-\alpha}{2}\right)\zeta(z)\zeta(z-\alpha)\frac{\tau^{(s-z)/2}}{z-s}\prod_{p \in \mathbb{P}}\left(1-\frac{1}{p^{z}} \right)\left(1-\frac{1}{p^{z-\a}} \right)dz.
	\end{align}
	Note that the functional equation for $\zeta(z)\zeta(z-\a)$ in variable $z$ can be easily obtained by multiplying the functional equations for $\zeta(z)$ and $\zeta(z-\a)$ as
	\begin{align*}
		\pi^{-z+\a/2}\Gamma\left( \frac{z}{2}\right) \Gamma\left( \frac{z-\a}{2}\right) \zeta(z)\zeta(z-\a)=\pi^{-(1-z) -\a/2}\Gamma\left(\frac{1-z}{2} \right)\Gamma\left(\frac{1-z+\a}{2} \right) \zeta(1-z)\zeta(1-z+\a).
	\end{align*}
	Employing the above functional relation in $\eqref{f_2}$, we find that 
	\begin{align}
		f_{2}=&\frac{\pi^{-1-\a/2}}{2\pi i}\int_{\s_1+i\i}^{\s_1-i\i}\pi^{z}\Gamma\left( \frac{1-z}{2}\right)\Gamma\left( \frac{1-z+\a}{2}\right)\times\nonumber\\
		&\zeta(1-z)\zeta(1-z+\alpha)	\frac{\tau^{(s-z)/2}}{z-s}\prod_{p \in \mathbb{P}}\left(1-\frac{1}{p^{z}} \right)\left(1-\frac{1}{p^{z-\a}} \right)dz.  \label{f2_1}
	\end{align}
	Our choice of $\s_1$ helps us to change the factor $\zeta(1-z)\zeta(1-z+\a)$ into a Dirichlet series,  
	\begin{align}
		\zeta(1-z)\zeta(1-z+\a) =\sum_{n=1}^{\i}\frac{\s_{\a}(n)}{n^{1-z+\a}}.\label{eq: dirser}
	\end{align}
	Thus using \eqref{eq: dirser}, and then interchanging the order of summation and integration in \eqref{f2_1}, we get
	\begin{align*}
		f_2=&\frac{\pi^{-1-\a/2}}{2\pi i}\sum_{n=1}^{\i}\frac{\s_\a(n)}{n^{\a+1}}\times \\
		&\int_{\s_1+i\i}^{\s_1-i\i}\pi^{z}\Gamma\left( \frac{1-z}{2}\right)\Gamma\left( \frac{1-z+\a}{2}\right)	\frac{\tau^{(s-z)/2}n^{z}}{z-s}\prod_{p \in \mathbb{P}}\left(1-\frac{1}{p^{z}} \right)\left(1-\frac{1}{p^{z-\a}} \right)dz.
	\end{align*}
	By applying the change of variable $ z$ by $s-z$, one obtains
	\begin{align*}
		f_2=&\frac{\pi^{s-1-\a/2}}{2\pi i}\sum_{n=1}^{\i}\frac{\s_\a(n)}{n^{\a+1-s}}\int_{\s-\s_1-i\i}^{\s-\s_1+i\i}\Gamma\left( \frac{1+z-s}{2}\right)\Gamma\left( \frac{1+z-s+\a}{2}\right)\\
		&\times\left( \frac{\sqrt{\tau}}{n\pi}\right)^{z}\prod_{p \in \mathbb{P}}\left(1-\frac{1}{p^{s-z}} \right)\left(1-\frac{1}{p^{s-z-\a}} \right)\frac{dz}{z}.
	\end{align*}	
	Invoking Lemma \ref{identity for mu}, with $z$ replaced by $ s-z$, 
	\begin{align}\label{f_2contour}
		f_2=&\frac{\pi^{s-1-\a/2}}{2\pi i}\sum_{k|\bar{\mathbb{P}}^{2}}\frac{(\mu \ast \mu N^\a)~(k)}{k^{s}}\sum_{n=1}^{\i}\frac{\s_\a(n)}{n^{\a+1-s}}\times\nonumber\\
		&\int_{\s-\s_1-i\i}^{\s-\s_1+i\i}\Gamma\left( \frac{1+z-s}{2}\right)\Gamma\left( \frac{1+z-s+\a}{2}\right)	\left( \frac{n\pi}{\sqrt{\tau}k}\right)^{-z}\dfrac{dz}{z}.
	\end{align}
	Now we want to apply Lemma \ref{Mellin 1.0}, with $ u=\frac{n\pi}{k\sqrt{\tau}}$,  $s\to 1-s$ and $\a \to -\a$. See that for $\textup{Re}(z)= \s- \s_1$, all the required conditions to employ Lemma \ref{Mellin 1.0} in \eqref{f_2contour} are satisfied. Thus we derive, for $\textup{Re}(\t)>0,$
	\begin{align}\label{beforeramanujanguinand}
		f_2=\frac{2}{\tau^{\a/4+(1-s)/2}}\sum_{k|\bar{\mathbb{P}}^{2}}\frac{(\mu \ast \mu N^\a)~(k)}{k^{1+\a/2}}\sum_{n=1}^{\i}\frac{\s_\a(n)}{n^{\a/2}}\int_{1}^{\i}t^{\a/4-(s+1)/2}K_{\a/2}\left( \frac{2n\pi}{k}\sqrt{\frac{t}{\tau}}\right)dt.
	\end{align}
	The asymptotic bound on $K_{\nu}(z)$ from Lemma \ref{asymptotic} will help us to interchange the summation and integration in $\eqref{beforeramanujanguinand}$, to obtain,
	\begin{align}
		f_2 =2\tau^{-\a/4-(1-s)/2}\int_{1}^{\i}t^{\a/4-(s+1)/2}\sum_{k|\bar{\mathbb{P}}^{2}}\frac{(\mu \ast \mu N^\a)~(k)}{k^{1+\a/2}}\sum_{n=1}^{\i}\frac{\s_\a(n)}{n^{\a/2}}K_{\a/2}\left( \frac{2n\pi}{k}\sqrt{\frac{t}{\tau}}\right)dt.
	\end{align}
	Invoking Theorem \ref{doublesumBessel} in the above equation, we have 
	\begin{align}\label{f_2main}
		f_2 =2\tau^{-\a/4-(1-s)/2}\int_{1}^{\i}t^{\a/4-(s+1)/2}\Bigg(\sqrt{\frac{\t}{t}}\sum\limits_{\mathop {n=1 {}}\limits_{(n,\bar{\mathbb{P}}^2)=1} }^{\i}&\frac{\s_\a(n)}{n^{\a/2}}K_{\a/2}\left( 2n\pi\sqrt{\frac{\t}{t}}\right)-\left(\mathcal{F}(\a)+\mathcal{F}(-\a) \right)\Bigg)\ dt,
	\end{align}
	where $\mathcal{F}(\a)$ is defined in \eqref{FG}.
	
	Upon using $\eqref{cauchy}, ~\text{and}~ \eqref{f_1main},~ \text{in}~\eqref{f_2main}$, and noting that $(n,\bar{\mathbb{P}}^2)=1$ if and only if $(n,\bar{\mathbb{P}})=1 $,  we derive
	\begin{align}
		2\tau^{s/2-\a/4}&\int_{1}^{\i}t^{-\a/4+s/2-1}\Bigg(\sum\limits_{\mathop {n=1 {}}\limits_{(n,\bar{\mathbb{P}})=1} }^{\i}\frac{\s_\a(n)}{n^{\a/2}}K_{\a/2}\left( 2n\pi \sqrt{\t t}\right)\nonumber\\
		& +\frac{t^{\a/2}}{t^s}\Bigg(\sum\limits_{\mathop {n=1 {}}\limits_{(n,\bar{\mathbb{P}})=1} }^{\i}\frac{\s_\a(n)}{n^{\a/2}}K_{\a/2}\left( 2n\pi\sqrt{\frac{\t}{t}}\right) 
		-\sqrt{\frac{t}{\t}}\left(\mathcal{F}_\a(t)+\mathcal{F}_{-\a}(t) \right)\Bigg)\Bigg)dt\nonumber\\ 
		&=\pi^{(\alpha-1)/2} \Gamma\left( \frac{1-\alpha}{2}\right)\zeta(1-\alpha)\frac{\tau^{(s-1)/2}}{1-s}\prod_{p \in \mathbb{P}}\left(1-\frac{1}{p} \right)\left(1-\frac{1}{p^{1-\a}} \right)\nonumber\\
		&+\pi^{-( 1+\a)/ 2}\Gamma\left( \frac{1+\a}{2}\right)\zeta(1+\a)\frac{\tau^{( s-1-\a) /2}}{1+\a-s}\prod_{p \in \mathbb{P}}\left(1-\frac{1}{p} \right)\left(1-\frac{1}{p^{1+\a}} \right)\nonumber\\
		&+\pi^{-s+\a/ 2}\Gamma\left( \frac{s}{2}\right)\Gamma\left( \frac{s-\a}{2}\right)\zeta(s)\zeta(s-\a)\prod_{p \in \mathbb{P}}\left(1-\frac{1}{p^{s}} \right)\left(1-\frac{1}{p^{s-\a}} \right).\label{aditi}
	\end{align}
From the left-hand side of the above equation \eqref{aditi}, we will evaluate the integration of the $n=1$ terms of both the series separately as:
\begin{align}\label{rajat}
\hspace{-1cm}\int_{1}^{\i}t^{-\a/4+s/2-1}K_{\a/2}\left( 2\pi \sqrt{\t t}\right)dt +\int_{1}^{\i}t^{\a/4-s/2-1}&K_{\a/2}\left( 2\pi \sqrt{\frac{\t}{ t}}\right)dt\nonumber\\
&=\int_{0}^{\infty}t^{-\a/4+s/2-1}K_{\a/2}\left( 2\pi \sqrt{\t t}\right)dt\nonumber\\
&=\frac{1}{2}(\pi \sqrt{\tau})^{-s+\alpha/2}\Gamma\left(\frac{s}{2}\right)\Gamma\left(\frac{s-\alpha}{2}\right),
\end{align}
provided $\textup{Re}(s)> \max\{0,\text{Re}(\alpha)\}$ and $\textup{Re}(\tau)>0$. In the last step we have used \cite[Formula~6.561.16, p.~676]{grn}, namely for $\textup{Re}(\nu+1\pm \mu)>0$ and $\textup{Re}(a)>0$, 
\begin{align*}
\int_{0}^{\infty}x^{\mu}K_{\nu}(ax)\ dx=\frac{2^{\mu-1}}{a^{\mu+1}}\Gamma\left(\frac{1+\mu+\nu}{2}\right)\Gamma\left(\frac{1+\mu-\nu}{2}\right).
\end{align*}

Now, dividing \eqref{aditi} by the factor $\displaystyle{\pi^{-s+\a/ 2}\Gamma\left( \frac{s}{2}\right)\Gamma\left( \frac{s-\a}{2}\right)}$ and simplifying the terms corresponding to $n=1$ by employing \eqref{rajat}, we derive \eqref{Maintheorem1eqn} for $\textup{Re}(s)> \max\{0,\text{Re}(\alpha)\}$. $s \neq 1, 1 +\a$, and $\textup{Re}(\tau)>0$.
This completes the proof. 
\end{proof}

Next, we will provide proof of Theorem \ref{Theoalpha0eqnintro}.

\subsection{Proofs of Theorem \textup{\ref{Theoalpha0eqnintro}} and Corollary \textup{\ref{coroalpha0>1}}}

\begin{proof}[Proof of Theorem \textup{\ref{Theoalpha0eqnintro}}]
Letting $\a \to 0$ in \eqref{Maintheorem1eqn}, we get 
\begin{align}\label{zetasquarefirststep}
	\zeta^2(s)\prod_{p \in \mathbb{P}}\left(1-\frac{1}{p^s} \right)^2&=1+\frac{\pi^{s}\tau^{s/2}}{\Gamma^2\left( \frac{s}{2}\right)}\Bigg[\lim_{\a \to 0}\left(\frac{\t^{-1/2}}{s-1}\mathcal{R}_{1-\a} +\frac{\t^{-(1+\a)/2}}{s-1-\a}
	\right)\nonumber \\
	&+2\int_{1}^{\i}t^{s/2-1}\Bigg\{\sum\limits_{\mathop {n=2 {}}\limits_{(n,\bar{\mathbb{P}})=1} }^{\i}d(n)K_{0}\left( 2n\pi \sqrt{\t t}\right) \nonumber\\
	& +t^{-s}\Bigg(\sum\limits_{\mathop {n=2 {}}\limits_{(n,\bar{\mathbb{P}})=1} }^{\i}d(n)K_{0}\left( 2n\pi\sqrt{\frac{\t}{t}}\right)-\sqrt{\frac{t}{\t}}\lim_{\a \to 0}\left(\mathcal{F}_\a(t)+\mathcal{F}_{-\a}(t) \right)\Bigg)\Bigg\}dt\Bigg].\nonumber \\
\end{align}
We first compute the limit,
\begin{align*}
	\lim_{\a \to 0}\left(\frac{\t^{-1/2}}{s-1}\mathcal{R}_{1-\a} +\frac{\t^{-(1+\a)/2}}{s-1-\a}\mathcal{R}_{1+\a} \right).
\end{align*}
To do so, we expand each factors involving $\alpha$ around $\a = 0$ 
to get,
\begin{align}\label{RlaphaBiggo}
	\mathcal{R}_{1-\a}&=\pi^{-(1-\a)/2}\Gamma\left(\frac{1-\a}{2}\right)\zeta(1-\alpha)\prod_{p \in  \mathbb{P}}\left(1-\frac{1}{p}\right)\left(1-\frac{1}{p^{1-\a}}\right)\nonumber\\
	&=\left( 1+
	\left(\frac{1}{2}\log \pi+ \log 2 +\frac{\gamma}{2}\right)\alpha
	+O(|\alpha|^2)\right)
	\bigg(-1/\alpha + \gamma + O(|\alpha|) \bigg) \nonumber\\
	&\qquad\qquad\qquad\qquad\qquad\qquad\qquad \times\bigg(\sum_{n|\bar{\mathbb{P}}}\frac{\mu(n)}{n}+\sum_{n|\bar{\mathbb{P}}}\frac{\mu(n)\log n}{n}\a +\mathcal{O}(|\a|^2)\bigg)\prod_{p \in \mathbb{P}}\left(1-\frac{1}{p} \right)\nonumber\\
	&=-\frac{1}{\a}\sum_{n|\bar{\mathbb{P}}}\frac{\mu(n)}{n}\prod_{p \in \mathbb{P}}\left(1-\frac{1}{p} \right)
	+\bigg(- \sum_{n|\bar{\mathbb{P}}}\frac{\mu(n)\log n}{n}+\left(\frac{\gamma}{2} -\log (\sqrt{4\pi})\right)\sum_{n|\bar{\mathbb{P}}}\frac{\mu(n)}{n}\bigg)
 \nonumber\\
	&\qquad\qquad\qquad\qquad\qquad\qquad\qquad\qquad\qquad\qquad\qquad\times \prod_{p \in \mathbb{P}}\left(1-\frac{1}{p} \right) +\mathcal{O}_{\bar{\mathbb{P}}}(|\a|).
\end{align}
In the first step above, we have used the Taylor series expansion
\begin{align*}
	\pi^{-(1-\alpha)/2}
	\Gamma\left(\frac{1-\alpha}{2}\right)
	=1+
	\left(\frac{1}{2}\log \pi+ \log 2 +\frac{\gamma}{2}\right)\alpha
	+O(|\alpha|^2),
\end{align*}
and the Laurent series expansion $\zeta(1 -\alpha) = -1/\alpha + \gamma + O(|\alpha|)$ around $\alpha = 0$. In order to expand the product over primes involving $\a$, we first write it in a
different form as 
\begin{align*}
	\prod_{p \in \mathbb{P}}\left(1-\frac{1}{p^{1-\a}} \right) =\sum_{n|\bar{\mathbb{P}}}\frac{\mu(n)}{n^{1-\a}}, 
\end{align*}
and then using the Laurent series expansion of $n^\a$ around $\a=0$, we obtain
\begin{align*}
	\prod_{p \in \mathbb{P}}\left(1-\frac{1}{p^{1-\a}} \right)= \left(\sum_{n|\bar{\mathbb{P}}}\frac{\mu(n)}{n}+\sum_{n|\bar{\mathbb{P}}}\frac{\mu(n)\log n}{n}\a +\mathcal{O}(|\a|^2)\right).
\end{align*}
Similarly in $\mathcal{R}_{1+\a}$, using the Taylor series of $\pi^{-(1+\alpha)/2}
\Gamma\left(\frac{1+\alpha}{2}\right)$, $\t^{-\a/2}$, $(s-1-\a)^{-1}$, and the Laurent series of $\zeta(1+\a)$  around $\alpha=0$, we have 
\begin{align}\label{R-laphaBiggo}
	\frac{\t^{-\a/2}}{s-1-\a}\mathcal{R}_{1+\a} &= \left( 1-\frac{\log \t}{2}\a +\mathcal{O}(|\a|^2)\right)\left(\frac{1}{s-1}+\frac{\a}{(s-1)^2}+\mathcal{O}(|\a|^2)\right)\nonumber\\
	&\qquad \qquad\times
	\left( 1-
	\left(\frac{1}{2}\log \pi+ \log 2 +\frac{\gamma}{2}\right)\alpha
	+O(|\alpha|^2)\right)\left(\frac{1}{\a}+\gamma+\mathcal{O}(|\a|) \right)\nonumber\\
	&\qquad \qquad\times\bigg(\sum_{n|\bar{\mathbb{P}}}\frac{\mu(n)}{n}-\sum_{n|\bar{\mathbb{P}}}\frac{\mu(n)\log n}{n}\a +\mathcal{O}(|\a|^2)\bigg)\prod_{p \in \mathbb{P}}\left(1-\frac{1}{p} \right)\nonumber \\
	&=\frac{1}{\a(s-1)}\sum_{n|\bar{\mathbb{P}}}\frac{\mu(n)}{n}\prod_{p \in \mathbb{P}}\left(1-\frac{1}{p} \right)+\frac{1}{(s-1)}\bigg(-\sum_{n|\bar{\mathbb{P}}}\frac{\mu(n)\log n}{n}\nonumber\\
	&\qquad\qquad+\left(\frac{1}{s-1} +\frac{\gamma}{2}-\log(\sqrt{4\pi \t})\right)\sum_{n|\bar{\mathbb{P}}}\frac{\mu(n)}{n}\bigg)\prod_{p \in \mathbb{P}}\left(1-\frac{1}{p} \right)+\mathcal{O}_{\bar{\mathbb{P}}}(|\a|).
\end{align}
Therefore, around $\alpha=0$, addition of \eqref{RlaphaBiggo} and \eqref{R-laphaBiggo} gives
\begin{align*}
	\frac{\t^{-1/2}}{s-1}\mathcal{R}_{1-\a} &+\frac{\t^{-(1+\a)/2}}{s-1-\a}\mathcal{R}_{1+\a}\\
	&= \frac{\tau^{-1/2}}{s-1}\prod_{p \in \mathbb{P}}\left(1-\frac{1}{p} \right)\left(-2\sum_{n|\bar{\mathbb{P}}}\frac{\mu(n)\log n}{n}+\left(\frac{1}{s-1}+\gamma-\log(4\pi \sqrt{\t}) \right)\sum_{n|\bar{\mathbb{P}}}\frac{\mu(n)}{n} \right) +\mathcal{O}_{\bar{\mathbb{P}}}(\a).
\end{align*}
Thus,
\begin{align}\label{limitalpha1}
	\lim_{\a \to 0}&\left(\frac{\t^{-1/2}}{s-1}\mathcal{R}_{1-\a} +\frac{\t^{-(1+\a)/2}}{s-1-\a}\mathcal{R}_{1+\a} \right) \nonumber \\
	&= \frac{\tau^{-1/2}}{s-1}\prod_{p \in \mathbb{P}}\left(1-\frac{1}{p} \right)\left(-2\sum_{n|\bar{\mathbb{P}}}\frac{\mu(n)\log n}{n}+\left(\frac{1}{s-1}+\gamma-\log(4\pi \sqrt{\t}) \right)\sum_{n|\bar{\mathbb{P}}}\frac{\mu(n)}{n} \right).
\end{align}
Similarly, we can find that, 
\begin{align}\label{limitalpha2}
	\lim_{\a \to 0}\left(\mathcal{F}_{\a}(t)+\mathcal{F}_{-\a}(t)\right)= \frac{1}{4}\prod_{p \in \mathbb{P}}\left(1-\frac{1}{p} \right)\left(- 2\sum_{n|\bar{\mathbb{P}}}\frac{\mu(n)\log n}{n}+\left(\gamma -\log\left(4\pi\sqrt{\frac{\t}{t}}\right) \right)\sum_{n|\bar{\mathbb{P}}}\frac{\mu(n)}{n} \right).
\end{align}
We conclude by employing \eqref{limitalpha1}, \eqref{limitalpha2} into \eqref{zetasquarefirststep}, and writing $\sum_{n|\bar{\mathbb{P}}}\frac{\mu(n)}{n}$ in the form of Dirichlet series $\prod_{p \in  \mathbb{P}}\left(1-\frac{1}{p} \right)$.
\end{proof}
Now, we turn our attention to prove Corollary \ref{coroalpha0>1}.

\begin{proof}[Proof of Corollary \textup{\ref{coroalpha0>1}}]
In Theorem \ref{Theoalpha0eqnintro}, we split the integration over $t$ on the right-hand side, and consider the integral, 
\begin{align}\label{inegralFmain}
	\int_{1}^{\i}t^{-s/2-1/2}&\left(-2\sum_{n|\bar{\mathbb{P}}}\frac{\mu(n)\log n}{n}+\left(\gamma -\log\left(4\pi\sqrt{\frac{\t}{t}}\right) \right)\prod_{p \in \mathbb{P}}\left(1-\frac{1}{p} \right) \right)dt.
\end{align}
Using integral identities, for  $\textup{Re}(s)>1$,
\begin{align}
	\int_{1}^{\i}t^{-s/2-1/2}dt &=\frac{2}{s-1}, \\
	\text{ and } \int_{1}^{\i}t^{-s/2-1/2}\log\left(4\pi \sqrt{\frac{\t}{t}} \right)dt&=-\frac{2}{(s-1)^2}+\frac{2}{(s-1)}\log\left(4\pi \sqrt{\t}\right), 
\end{align}
in \eqref{inegralFmain}, we obtain, 
\begin{align*}
	\int_{1}^{\i}t^{-s/2-1/2}&\left(-2\sum_{n|\bar{\mathbb{P}}}\frac{\mu(n)\log n}{n}+\left(\gamma -\log\left(4\pi\sqrt{\frac{\t}{t}}\right) \right)\prod_{p \in \mathbb{P}}\left(1-\frac{1}{p} \right) \right) dt\nonumber\\
	&=-\frac{4}{(s-1)}\sum_{n|\bar{\mathbb{P}}}\frac{\mu(n)\log n}{n}+2\left(\frac{\gamma}{(s-1)}+\frac{1}{(s-1)^2}-\frac{1}{(s-1)}\log\left(4\pi \sqrt{\t}\right)\right)\prod_{p \in \mathbb{P}}\left(1-\frac{1}{p} \right).
\end{align*}
Thus, for $\textup{Re}(s)>1$, and $\textup{Re}(\t)>0$,
\begin{align}\label{Finalintegral}
	-&\frac{1}{2\sqrt{\t}}\prod_{p \in \mathbb{P}}\left(1-\frac{1}{p} \right)\int_{1}^{\i}t^{-s/2-1/2}\left(-2\sum_{n|\bar{\mathbb{P}}}\frac{\mu(n)\log n}{n}+\left(\gamma -\log\left(4\pi\sqrt{\frac{\t}{t}}\right) \right)\prod_{p \in \mathbb{P}}\left(1-\frac{1}{p} \right) \right)dt\nonumber\\
	&=-\frac{1}{(s-1)\sqrt{\t}}\prod_{p \in \mathbb{P}}\left(1-\frac{1}{p} \right)\left(-2\sum_{n|\bar{\mathbb{P}}}\frac{\mu(n)}{n}\log n +\left(\frac{1}{s-1}+\gamma -\log(4\pi \sqrt{\t})\right)\prod_{p \in\mathbb{P}}\left(1-\frac{1}{p} \right)\right).
\end{align}
From $\eqref{Finalintegral}$ we observe that the second and the last term on the right-hand side of \eqref{coroalpha0eqnintro} beautifully cancels out each other. Finally, the expression in \eqref{coroalpha0eqnintro} vastly simplifies to the required equation \eqref{eq},
for $\textup{Re}(s)>1$, and $\textup{Re}(\t)>0.$
\end{proof}

\section{Acknowledgment}
We are grateful to the reviewer who gave essential comments and feedback to improve our article. The authors would also like to thank Prof. Atul Dixit and Prof. Akshaa Vatwani for the insightful discussions and the critical first review of this manuscript.

%
%
%

\end{document}